\documentclass[a4paper,11pt]{article}
\usepackage[margin=1in]{geometry}
\makeatletter
\usepackage{palatino}
\newcommand*{\rom}[1]{\expandafter\@slowromancap\romannumeral #1@}
\makeatother
\usepackage{hyperref}
\hypersetup{ colorlinks=true, citecolor=black, linkcolor=black, urlcolor=black} 

\def\a{\alpha}

\def\k{\kappa}
\def\l{\lambda}
\def\m{\mu}

\usepackage{amsfonts}
\usepackage{amsmath}
\usepackage{amssymb}
\usepackage{lmodern}
\usepackage[T1]{fontenc}
\usepackage{pb-diagram}
\usepackage[cmtip,arrow]{xy}
\usepackage{pb-diagram,pb-xy}
\usepackage{pgf}
\usepackage{epigraph}
\usepackage{graphicx}

\usepackage{tikz}
\usetikzlibrary{arrows,decorations.pathmorphing,backgrounds,positioning,fit,calc}
\usepackage{tikz-cd}
\usepackage{float}

\usepackage{pb-diagram}
\usepackage[cmtip,arrow]{xy}
\usepackage{pb-diagram,pb-xy}
\usepackage{pgf}
\usepackage{epigraph}
\usepackage{graphicx}
\usepackage{amsmath}
\usepackage{amssymb}
\usepackage{tikz}
\usepackage{tikz-cd}
\usepackage{enumerate}
\usepackage{amsfonts}
\usepackage{comment}
\usepackage{amsthm}
\usepackage{amssymb}
\usepackage{verbatim}
\usepackage{mathtools}
\usepackage{titlesec}
\usepackage{chngcntr}
\usepackage{mathrsfs}

\titleformat{\subsubsection}[runin]{}{\thesubsubsection}{}{}[]

\titleformat{\subsubsection}[runin]
{\normalfont\bfseries}
{\thesubsubsection}{0.5em}{}

\titleformat{\section}[block]{}{\thesubsection}{}{}[]

\titleformat{\section}[block]
{\bfseries \huge}
{\thesection}{0.5em}{}

\titleformat{\subsection}[block]{}{\thesubsection}{}{}[]

\titleformat{\subsection}[block]
{\bfseries \Large}
{\thesection}{0.5em}{}

\counterwithin*{subsubsection}{subsection}
\counterwithin*{subsubsection}{section}
\newtheorem{thm}{Theorem}[section]

\newtheorem{defn}[thm]{Definition}
\newtheorem{prop}[thm]{Proposition}
\newtheorem{cor}[thm]{Corollary}
\newtheorem{lem}[thm]{Lemma}
\newtheorem{conj}[thm]{Conjecture}
\newtheorem{example}[thm]{Example}
\newtheorem{rmk}[thm]{Remark}
\newtheorem{assertion}{ \color{red} Assertion}

\newtheorem{idea}{Idea}[section]
\newtheorem{aim}[thm]{Aim}
\newtheorem{qn}[thm]{Question}
\newtheorem{cl}[subsubsection]{Claim}
\newtheorem{ntt}[subsubsection]{Notation}
\newtheorem{nt}[subsubsection]{Note}

\newtheorem{slg}{Slogan}

\newenvironment{reference}{\paragraph{Reference:}}{\hfill$\square$}

\newcommand{\brf}{\begin{reference}}
\newcommand{\erf}{\end{reference}}
\newcommand{\bnt}{\begin{nt} \normalfont}
\newcommand{\ent}{\end{nt}}
\newcommand{\bntt}{\begin{ntt}}
\newcommand{\entt}{\end{ntt}}
\newcommand{\bcl}{\begin{cl} \normalfont }
\newcommand{\ecl}{\end{cl}}
\newcommand{\bqn}{\begin{qn} \normalfont  }
\newcommand{\eqn}{ \end{qn}}
\newcommand{\bid}{\begin{idea}}
\newcommand{\eid}{\end{idea}}

\newcommand{\bas}{\begin{assertion} \normalfont \begin{bf} }
\newcommand{\eas}{\end{bf} \end{assertion}}
\newcommand{\bcr}{\begin{cor}}
\newcommand{\ecr}{\end{cor}}
\newcommand{\bex}{\begin{example} \normalfont}
\newcommand{\eex}{\end{example}}
\newcommand{\blm}{\begin{lem}}
\newcommand{\elm}{\end{lem}}
\newcommand{\bthm}{\begin{thm}}
\newcommand{\ethm}{\end{thm}}
\newcommand{\bcd}{\begin{tikzcd}}
\newcommand{\ecd}{\end{tikzcd}}
\newcommand{\bdf}{\begin{defn} \normalfont}
\newcommand{\edf}{\end{defn}}
\newcommand{\bpp}{\begin{prop}}
\newcommand{\bam}{\begin{aim}}
\newcommand{\eam}{\end{aim}}

\newcommand{\rar}{\rightarrow}
\newcommand{\brm}{\begin{rmk} \normalfont }
\newcommand{\erm}{\end{rmk}}
\newcommand{\epp}{\end{prop}}
\newcommand{\bpf}{\begin{proof}}

\newcommand{\hrar}{\hookrightarrow}

\newcommand{\epf}{\end{proof}}

\newcommand{\kom}{$\mathcal{\kom}$}
\newcommand{\CC}{\mathbb{C}}

\newcommand{\rd}{\textcolor{red}}

\newcommand{\pno}{\mathbb{P}^{1}}

\newcommand{\overbar}[1]{\mkern 1.5mu\overline{\mkern-1.5mu#1\mkern-1.5mu}\mkern 1.5mu}

\newcommand{\bnu}{\begin{enumerate}}
\newcommand{\enu}{\end{enumerate}}

\newcommand{\PP}{\mathbb{P}}
\newcommand{\ZZ}{\mathbb{Z}}
\newcommand{\NN}{\mathbb{N}}
\newcommand{\QQ}{\mathbb{Q}}
\newcommand{\Gm}{\mathbb{G}_{m}}
\newcommand{\Ga}{\mathbb{G}_{a}}
\newcommand{\bdoc}{\begin{document}}
\newcommand{\edoc}{\end{document}}
\newcommand{\bpm}{\begin{pmatrix}}
\newcommand{\epm}{\end{pmatrix}}
\newcommand{\eps}{\varepsilon}
\newcommand{\bfct}{\begin{fact}}
\newcommand{\efct}{\end{fact}}
\newcommand{\bslg}{\begin{slg}}
\newcommand{\eslg}{\end{slg}}
\newcommand{\Uh}{\hat{U}}

\newcommand{\hU}{\hat{U}}
\newcommand{\git}{\mathbin{
		\mathchoice{/\mkern-6mu/}
		{/\mkern-6mu/}
		{/\mkern-5mu/}
		{/\mkern-5mu/}}}
\newcommand{\mc}{\mathscr}
\newcommand{\rqq}{\rd{[??]}}
\newcommand{\re}{\rd{[...]}}
\DeclareMathOperator{\Mat}{Mat}
\newcommand{\diagentry}[1]{\mathmakebox[1.8em]{#1}}
\newcommand{\xddots}{%
	\raise 4pt \hbox {.}
	\mkern 6mu
	\raise 1pt \hbox {.}
	\mkern 6mu
	\raise -2pt \hbox {.}
}
\newcommand{\bgs}{\begin{gss}}
\newcommand{\egs}{\end{gss}}

\DeclareMathOperator{\GL}{GL}
\DeclareMathOperator{\SL}{SL}
\DeclareMathOperator{\PGL}{PGL}
\DeclareMathOperator{\SU}{SU}
\DeclareMathOperator{\Sp}{Sp}
\DeclareMathOperator{\SO}{SO}
\DeclareMathOperator{\U}{U}
\DeclareMathOperator{\proj}{Proj}
\DeclareMathOperator{\spec}{Spec}
\DeclareMathOperator{\Hom}{Hom}
\DeclareMathOperator{\dom}{dom}
\DeclareMathOperator{\diag}{diag}
\DeclareMathOperator{\trace}{trace}
\DeclareMathOperator{\sym}{Sym}
\DeclareMathOperator{\conv}{Conv}
\DeclareMathOperator{\stab}{Stab}
\DeclareMathOperator{\im}{\text{Im}}
\DeclareMathOperator{\wt}{wt}
\DeclareMathOperator{\mult}{mult}
\DeclareMathOperator{\id}{id}
\DeclareMathOperator{\pr}{pr}
\DeclareMathOperator{\ev}{ev}
\DeclareMathOperator{\res}{res}
\DeclareMathOperator{\spn}{Span}
\DeclareMathOperator{\Aut}{Aut}
\DeclareMathOperator{\Lie}{Lie}
\DeclareMathOperator{\Pic}{Pic}
\DeclareMathOperator{\Quot}{Quot}
\DeclareMathOperator{\End}{End}
\DeclareMathOperator{\Ext}{Ext}

\newcommand{\Zb}{Z_\beta}
\newcommand{\Yb}{Y_\beta}
\newcommand{\Sb}{S_\beta}
\newcommand{\Ybc}{Y_\beta^{c}}
\newcommand{\Zbc}{Z_\beta^{c}}
\newcommand{\Pb}{P_\beta}

\DeclarePairedDelimiter\ceil{\lceil}{\rceil}
\DeclarePairedDelimiter\floor{\lfloor}{\rfloor}
\newcommand{\fn}{\floor*{\frac{n}{2}}}

\renewcommand{\Gamma}{\varGamma}
\renewcommand{\phi}{\varphi}
\renewcommand{\vec}[1]{\underline{#1}}
\renewcommand{\geq}{\geqslant}
\renewcommand{\leq}{\leqslant}

\mathchardef\ordinarycolon\mathcode`\:
\mathcode`\:=\string"8000
\begingroup \catcode`\:=\active
  \gdef:{\mathrel{\mathop\ordinarycolon}}
\endgroup
\renewcommand{\epigraphflush}{flushright}
\newcommand{\bcon}{\begin{conj}}
\newcommand{\econ}{\end{conj}}
\newcommand{\Hh}{\hat{H}}

\newcommand{\HN}{Harder-Narasimhan }
\newcommand{\sh}{\mc{F}}

\titleformat{\subsubsection}[runin]{}{\thesubsubsection}{}{}[]

\titleformat{\subsubsection}[runin]
{\normalfont\bfseries}
{\thesubsubsection}{0.5em}{}

\titleformat{\section}[block]{}{\thesubsection}{}{}[]

\titleformat{\section}[block]
{\bfseries \LARGE}
{\thesection}{0.5em}{}

\titleformat{\subsection}[block]{}{\thesubsection}{}{}[]

\titleformat{\subsection}[block]
{\bfseries \Large}
{\thesubsection}{0.5em}{}

\newcommand{\?}{\rd{???}}

\newcommand{\sss}{ \subsubsection{} }

\begin{document}
\title{ {\Huge Variation of Non-Reductive Geometric Invariant Theory}}  \author{Gergely B\'{e}rczi, Joshua Jackson, Frances Kirwan} 
\maketitle


Variation of Geometric Invariant Theory (VGIT) \cite{Dolgachev1998, Thaddeus1996} studies the structure of the dependence of a GIT quotient on the choice of linearisation. This structure, and the concomitant wall-crossing picture relating the different quotients when a reductive group acts linearly on a projective variety (with respect to an ample linearisation), has long been a hallmark of classical GIT, and has found many diverse applications. In this article we show that under the conditions of results given in \cite{Berczi2016a} for non-reductive linear algebraic group actions on projective varieties, this structure persists in non-reductive GIT.
\\~\\
In $\S$ \ref{review} we review the key points of the classical theory, when a reductive linear algebraic group $G$ acts on a projective variety $X$. Mumford's GIT associates  to any linearisation $\mathcal{L}$ of this action with respect to an ample line bundle $L$  a  notion of a quotient $X/\!/_{\mathcal{L}} G$ (where $X/\!/_{\mathcal{L}} G$ is a projective variety), and the variation of GIT results of Thaddeus \cite{Thaddeus1996} and Dolgachev \& Hu  \cite{Dolgachev1998} describe the dependence of $X/\!/_{\mathcal{L}} G$
on the linearisation $\mathcal{L}$. In $\S$ \ref{NRGIT} we give a brief exposition of how some of the difficulties characteristic of non-reductive GIT may be solved 
for actions of  linear algebraic groups with graded unipotent radicals as in  \cite{Berczi2016a} and \cite{Berczi2016}. The next two sections are dedicated to showing that these results allow us to recover a variation picture which is very similar to the classical one, in the case where  \lq semistability coincides with stability for the unipotent radical\rq , in a sense that we will specify later. The major difference from usual VGIT is the presence, a priori, of an additional parameter: a choice of a suitable \text{ }1-parameter subgroup of the group in question.
\\~\\
We then discuss what happens without the simplifying assumption that semistability coincides with stability for the unipotent radical. Essentially the same description can be made to work, with some slight modifications. Finally, we discuss some illustrative examples, and indicate some potential applications of our results.

\bigskip

 \begin{center} \textbf{Conventions, Notation and Standing Assumptions} \end{center} ~\\
Our ground field is $\CC$ (or more generally an algebraically closed field of characteristic 0). Throughout, unless otherwise specified, let $X$ be a irreducible projective variety carrying an ample line bundle $L$. When $X$ is acted on by a linear algebraic group, we say that the action is linear, or linearised, if we have chosen a linearisation: that is, a lift of the action on $X$ to the total space of $L$. We will denote linearisations themselves with cursive scripts, i.e. $\mc{L}$ is a linearisation with underlying bundle $L$. Unless otherwise stated, we will denote by $G$ a reductive linear algebraic group, and by $H$ an arbitrary linear algebraic group, which may be non-reductive. We fix a Levi factor $R$ of $H$, so that $H = U \rtimes R$ where $U$ is the unipotent radical of $H$. We will abbreviate \lq one-parameter subgroup\rq\text{ }to \lq 1PS\rq\text{ }in places. When given without qualification, point means closed point and sheaf means coherent sheaf.

\section{Reductive GIT}\label{review} \label{GIT}

When a reductive linear algebraic group $G$ acts on a projective variety $X$, Mumford's Geometric Invariant Theory allows us to construct  quotients, in the following sense. We need to make a choice of linearisation $\mc{L}$; that is, an ample line bundle $L$ on $X$, and a  lift of the $G$ action to $L$. Then the algebra   $\bigoplus\limits_{n=0}^\infty H^0(X, L^{\otimes n})^G$ of invariant sections is a finitely generated graded algebra; we call the associated projective variety $$X \git_{\mc{L}} G \coloneqq \mathrm{Proj}\bigoplus\limits_{n=0}^\infty H^0(X, L^{\otimes n})^G$$  the GIT quotient of $X$ by $G$ with respect to $\mc{L}$. Further, there are open subvarieties  $$X^{s}(\mc{L}) \subseteq X^{ss}(\mc{L}) \subseteq X,$$ respectively called the stable and semistable loci, such that the inclusion of the invariants into the algebra of sections induces a good quotient $X^{ss}(\mc{L}) \rar X\git_{\mc{L}}G$, which restricts to a geometric quotient $X^{s}(\mc{L}) \rar X^{s}(\mc{L})/G$ onto an open subvariety of the GIT quotient \cite{Mumford1994}. 

\brm When considering how the choice of linearisation affects the GIT quotient  it is important to note that, since good quotients are in particular categorical, if two linearisations yield the same semistable locus then their associated quotients are canonically isomorphic. \erm

Thus classical GIT provides us with $G$-invariant open subvarieties of the variety  $X$ and gives us the best possible quotients of those open subvarieties by the $G$-action. A crucial feature of GIT is that those quotients are not too hard to compute:  in order to determine the GIT quotient and its geometric restriction we need only determine the stable and semistable loci, and the following tells us that these are determined by essentially combinatorial data.

\bthm \cite{Mumford1994} \label{HMcrit} \textbf{ (Hilbert-Mumford Criterion)}  Let $T \subseteq G$ be a maximal torus, and choose coordinates that diagonalise the action of $T$. 

\bnu
\item Stability and semistability for the $T$-action on $X$ are determined purely combinatorially as follows. Take $x \in X$ and let $\Delta_x \subset \mathfrak{t}^*$ be the convex hull of the weights of the $T$-action corresponding to nonzero coordinates of $x$. Then $$x \in X^{ss, T}(\mc{L}) \text{ iff } 0 \in  \overbar{\Delta_x}; $$ $$x \in X^{s, T}(\mc{L}) \text{ iff } 0 \in \Delta_x. $$
\item A point $x \in X$ is semistable (resp. stable) for the $G$ action with respect to $\mc{L}$ iff every point of $G$-orbit is semistable (resp. stable) for $T$. That is, we have $$X^{(s)s, G}(\mc{L}) = \bigcap\limits_{g \in G} g \cdot X^{(s)s, T}(\mc{L}).$$

\enu
 \ethm

The study of how the GIT quotient depends on the linearisation, known as \emph{Variation of GIT} or VGIT, was initiated by Thaddeus in \cite{Thaddeus1996} and by Dolgachev and Hu in \cite{Dolgachev1998}.
\\~\\
Let $\text{Pic}^G(X)$ denote the set of $G$-linearised line bundles on $X$ up to the natural notion of isomorphism. The tensor product of two $G$-linearised line bundles and the dual of a $G$-linearised line bundle both have canonical $G$-linearisations, and hence this set has an abelian group structure. The following definition comes from \cite{Thaddeus1996}.

\bdf We say that two linearisations, $\mc{L}_1$, $\mc{L}_2$, are $G$-algebraically equivalent if there is a connected variety $S$ with points $s_1$,$s_2 \in S$, and a linearisation $\mc{L}_S$ of the $G$-action on $S \times X$ induced from the second factor, such that $\mc{L}_S \mathord \mid_{s_i} \cong \mc{L}_i$ for $i = 1,2.$ \edf

Thaddeus proved in \cite{Thaddeus1996} that the semistable locus in $X$, and hence the GIT quotient of $X$ by $G$, is invariant under the equivalence relation of $G$-algebraic equivalence, so we may pass to the quotient of $\text{Pic}^G(X)$ by this equivalence relation, denoted $\text{NS}^G(X)$. It follows from the Hilbert-Mumford criterion that the GIT quotient is unchanged if one replaces a given $G$-linearisation $\mc{L}$, of the line bundle $L$, with the induced linearisation $\mc{L}^n$ on $L^{\otimes n}$ for $n>0$. So we may consider the rational vector space $NS^G(X)\otimes \QQ$, which is shown to be finite dimensional. This is the space of rational linearisations. 

\bdf Let $A^G_{\QQ}$ be the cone in $NS^G(X)\otimes \QQ$ which consists of rational linearisation classes for which the underlying algebraic equivalence class of bundle is ample. Denote by $C^{G} \subset A^G_{\QQ}$ the set of $G$-effective rational ample linearisations; that is, those linearisations in $A^G_{\QQ}$ for which there is a semistable point in $X$. \edf

The goal of VGIT is to understand how the quotient changes as we move around in this parameter space. This is described in \cite{Dolgachev1998} using a so-called wall-and-chamber structure on $C^G$; that is, a subdivision of $C^G$ into convex chambers by homogeneous walls which are locally polyhedral of codimension one in $C^G$. If we assume that generic points of $X$ have 0-dimensional stabilisers, then a linearisation $\mc{L}$ lies on some {wall} if  it has a strictly semistable point: that is, there is a point of $X$ that is semistable but not stable with respect to $\mc{L}$. A connected component of the set of linearisations not lying on any wall is a \emph{chamber}.  A \emph{cell} is a connected component of the set of points of a given wall $W$ that do not lie on any other wall except those walls containing the whole of $W$. 

The main features of the dependence may be summarised as follows.

\bthm \cite{Dolgachev1998} Assume that there exists a linearisation such that semistability equals stability. Then:
\bnu  
\item The wall-and-chamber structure of $C^G$ is finite, in the sense that there are only finitely many chambers, walls and cells, and polyhedral, in the sense that the closure of a chamber is a rational polyhedral cone in the interior of $C^G$.
\item Two linearisations in the same chamber give the same semistable locus, and hence the same GIT quotient,
while the induced (rational) ample line bundle on the quotient varies affinely in the chamber. The set of linearisations giving a particular semistable locus is either a chamber or a union of cells in the same wall.

\item As a wall is crossed, the quotients for the wall and for the chambers on either side are related by an explicit birational transformation, called a Thaddeus flip. 
\enu  \ethm 

\brm The details of the birational transformation relating the quotients on either side of a wall will not concern us here. The definition may be found in \cite{Thaddeus1996}, where it is referred to simply as a flip. \erm

\bex A simple example of this structure is that of $\Gm$ acting linearly on $\mathbb{P}^n$. Then $C^G$ is the half-plane in $\QQ^2$, with the direction perpendicular to $\mathrm{NS}\otimes\QQ \cong \QQ$ corresponding to a choice of rational character. Choosing coordinates $x_0,..,x_n$ to diagonalise the action, and letting $\alpha_i \in \ZZ$ be the weight of $x_i$, we see that the chambers are bounded by rays passing through $(1,\alpha_i)$.     \eex

\section{Non-reductive GIT: the $\Uh$ Theorems} \label{NRGIT}

It is well known that many of the good properties of Mumford's GIT fail in general in the non-reductive setting. These include finite generation of the algebra of invariant sections, surjectivity of the quotient map, and the Hilbert-Mumford criteria. Recently, however, it has been shown that these properties hold for non-reductive linear algebraic group actions under special circumstances which are satisfied in a wide range of interesting applications. Here we present these results, collectively referred to as the various versions of the $\Uh$ theorem. We first introduce some notation.

\bdf Let $H$ be an arbitrary linear algebraic group. Let $U$ be its unipotent radical. We say that $H$ has \emph{graded unipotent radical} if there exists a one-parameter group of automorphisms $\Gm \leq \Aut H$ which acts with strictly positive weights on $\Lie U$, such that the induced action on $R=H/U$ is trivial. Such a $\Gm$ is called \emph{admissible}, and we refer to any given fixed choice as the \emph{distinguished} $\Gm$. We call an admissible $\Gm$ \emph{internal} if it comes from a subgroup $\Gm \leq H$, with the automorphism being conjugation, and \emph{external} otherwise. For the rest of this section we will assume that we have an internal admissible $\Gm$ (but see Remark \ref{extrm} below).
We will write $\Uh = U \rtimes \Gm$, or, for example, $\Uh_\lambda$ if we wish to emphasise a particular choice of homomorphism $\lambda : \Gm \rar \Aut(H)$. 

Up to conjugacy, we may assume that an internal admissible $\Gm$ is chosen from the Lie algebra $\mathfrak{z}$ of the centre of $R$, inside a maximal torus whose Lie algebra we denote by $\mathfrak{t}$. Then, choosing an invariant inner product, the admissible $\Gm$'s are precisely those whose pairings with all weights of the adjoint representation of $ \mathfrak{t}$ on $\Lie U$ are strictly positive. This defines a cone $\mathfrak{C} \subseteq \mathfrak{z}^\vee$, which we call the \emph{admissible cone}. \edf

\brm We emphasise that so far all of this is intrinsic to the group; it depends neither on $X$, nor on the action. \erm

\bdf Now suppose $H$ acts on $X$. We say that a linearisation $\mc{L}$ of an action of $H$ on $X$ with respect to an ample line bundle $L$ is \emph{adapted} if the lowest weight of the action of the distinguished $\Gm$ on $X$ is strictly negative, and the rest are strictly positive. Let $V_{min}^{\Gm}$ be the minimal weight space of the linear representation of the distinguished $\Gm$. Let $X_{\min}^{\Gm} = X \cap \mathbb{P}(V_{\min}^{\Gm})$, and let $X^{\text{o},{\Gm}}_{\min}$ be the subvariety consisting of those $x \in X$ with at least one coordinate associated to the minimal weight space being non-zero: that is, the basin of attraction of $X_{\min}^{\Gm}$ under the distinguished $\Gm$. We will observe later that if the unipotent radical $U$ of $H$ is non-trivial and $X$ is irreducible then $X_{\min}^{\Gm}$ and $X^{\text{o},{\Gm}}_{\min}$ are in fact independent of the choice of admissible $\Gm$ and so we will use the notation $X_{\min}$ and $X^{\text{o}}_{\min}$ instead. \edf

\brm If a given linearisation is not adapted, we can always twist the linearisation by a rational character to make it so, provided appropriate characters exist. \erm

\brm There is a technical point to be explained here, about the choice of linearisation. For the proofs of finite generation of invariants given in \cite{Berczi2016,Berczi2016a} to work, we must twist the linearisation by a rational character so that it is not merely adapted, but within some sufficiently small $\eps >0$ of the lower wall of the adapted chamber. This can make the statement of the results that follow somewhat unwieldy, so we adopt the following convention used in \cite{Berczi2016,Berczi2016a}: we say that a linearisation is \emph{well-adapted} if it is within distance $\eps >0$ of the relevant wall. More precisely, when in what follows we say that a property holds for a well-adapted linearisation, we mean that for any adapted linearisation there exists an $\eps >0$ such that the relevant property holds after twisting by a rational character so that the adapted linearisation lies within an open $\eps$-neighbourhood of the lower wall of the adapted chamber. This will only be relevant for finite generation of invariants; the stable and semistable loci and the associated quotients require simply that the linearisation should be adapted.
\erm

We are now ready to state the relevant versions of the $\Uh$ theorems.

\bthm \cite{Berczi2016a} \label{Uh} \textbf{($\Uh$-theorem when semistability coincides with stability for the unipotent radical)} Let $H$ act on a projective variety $X$, let $\Uh$ be formed as above by an admissible internal $\Gm$, and let $\mc{L}$ be a well-adapted linearisation with respect to a very ample line bundle $L$.  Assume that $$(*) \stab_U(x) = \{e\} \text{ for all } x \in X_{min}. $$ Then 

\bnu 
\item The $\Uh$-invariants are finitely generated, and the inclusion of the $\Uh$-invariant algebra induces a projective geometric quotient of an open subvariety $ X^{s, \Uh} =  X^{ss, \Uh}$ of $X$ (the $\Uh$-(semi)stable locus): $$ X^{s, \Uh} \rar X \git_\mc{L} \Uh. $$
\item Consequently the $H$-invariants are finitely generated, and the inclusion of the $H$-invariant subalgebra induces a good quotient of an open subvariety of $X$ (the $H$-semistable locus): $$X^{ss, H} \rar X \git_{\mc{L}} H,$$
where $X \git_{\mc{L}} H$ is the GIT quotient of $X \git_{\mc{L}} \Uh$ by the induced action of the reductive group $H/\Uh \cong R/\Gm $ with respect to the induced linearisation. The good quotient  $X^{ss, H} \rar X \git_{\mc{L}} H$ restricts to a geometric quotient $X^{s, H} \rar X^{s, H}/H \subseteq X \git_{\mc{L}} H$ of the $H$-stable locus $X^{s, H}$.
\item  \textbf{(Non-Reductive Hilbert-Mumford Criterion)} \label{NRHM} $x \in X$ is (semi)stable for $H$ if and only if it is (semi)stable for every $1$PS of $H$.  That is, if $T$ denotes any maximal torus of $H$, we have  $$X^{(s)s,H}(\mc{L}) = \bigcap\limits_{h \in H} hX^{(s)s,T}(\mc{L}).$$  
\enu
\ethm

\brm The condition $(*)$ holds iff we have $\stab_{U}(x) = \{e\}$ for all $x \in X^0_{\min}$. \erm

If the relevant stabilisers are not trivial for all points in the minimal weight space, we can proceed by a process analogous to, and indeed a generalisation of, the partial desingularisation process of \cite{kirw1985} for the reductive setting, when  a linearisation has strictly semistable points.

\bthm \cite{Berczi2016a} \label{Uh2}
 \textbf{($\Uh$-theorem giving projective completions)} With notation as in Theorem \ref{Uh}, suppose now that $(*)$ may fail, but we have instead the weaker condition that $$(**) \stab_U(x) = \{e\} \text{ for generic } x \in X_{\text{min}}.$$ Then

 \bnu 
 \item There exists a sequence of blow-ups of $X$ along $H$-invariant projective subvarieties resulting in a projective variety $\widehat{X}$ (with blow-down map $\widehat{\pi}: \widehat{X} \to X$) for which the conditions of Theorem \ref{Uh} hold for a suitable linearisation which can be taken to be an arbitrarily   small perturbation $\widehat{\mc{L}}$ of $\widehat{\pi}^*(\mc{L})$.
 \item There exists a further sequence of blow-ups along $H$-invariant projective subvarieties, resulting in a projective variety $\tilde{X}$ (with blow-down map $\tilde{\pi}: \tilde{X} \to X$)  such that the conditions of Theorem \ref{Uh} still hold for a suitable linearisation, and such that the quotient  given by that theorem is a geometric quotient of an open subvariety $\tilde{X}^{s,H}$ of $\tilde{X}$.
\enu
\ethm

\brm
This gives us an $H$-invariant open subvariety $X^{s,H}$ of $X$ with a geometric quotient by $H$ which is an open subvariety of the projective variety $\widehat{X} \git_{\widehat{\mc{L}}} H$ (and also of $\tilde{X} \git_{\tilde{\mc{L}}} H$). Here $X^{s,H}$ is the image under $\widehat{\pi}$ of the intersection of $\widehat{X}^{s,H}$ with the complement of the exceptional divisor in $\widehat{X}$.

The centres of the blow-ups used to obtain $\widehat{X}$ from $X$ are determined by the dimensions of the stabilisers in $U$ of the limits $\lim_{t \to 0} x$ for $x \in X^{0}_{\min}$. Another characterisation of $X^{s,\hU}$ is as 
$$X^{s,\hU} = \{ x \in \pi (\widetilde{X}^{s,\hU}) \mid \dim \stab_U(\lim_{t \to 0} x) = 0 \}. $$
\erm

There are applications for which all points of $X$ have non-trivial stabilisers in $U$, but we would still like to be able to perform quotients. There is a procedure that  allows this, yielding a further version of the $\Uh$ theorem.

\bthm \cite{Berczi2016a} \textbf{($\Uh$-theorem with positive-dimensional stabilisers in $U$)}
\label{Uh3}
Suppose that for the derived series  $$ U \geqslant U^{(1)} \geqslant ... \geqslant U^{(s)}  \geqslant \{e\}$$ of $U$, we have $$(***) \,\,\, \forall j \in \{1, \ldots, s\} \text{ }\exists d_j \in \NN \,\,\, \text{such that } \dim \stab_{U^{(j)}}(x) = d_j \text{ for all } x \in X^{\text{o}}_{\min}.$$ Then the conclusions of Theorem \ref{Uh} hold. 
\ethm

\brm In  fact the theorem may be applied with any series $U \geqslant U^{(1)} \geqslant ... \geqslant U^{(s)}  \geqslant \{e\}$ which is normalised by $H$ and has successive quotients $U^{(j)}/U^{(j+1)}$ abelian, provided that $(***)$ holds. \erm 

Finally, we may not have constant-dimensional stabilisers across $X^{0}_{\min}$ for any such series. In this case there is a final version of the $\Uh$ theorem.

\bthm \cite{Berczi2016a} \label{Uh4} \textbf{($\Uh$-theorem with positive-dimensional stabilisers in $U$, giving projective completions)}

There is a sequence of iterated blowups of $X$ along $H$-invariant closed subvarieties resulting in  $\pi: \widehat{X} \to X$ where $\widehat{X}$ is a projective variety with an induced linear action of $H$ for which $(***)$ holds for each $U^{(j)}$ in turn, and hence for which the conclusions of Theorem \ref{Uh3} hold.
\ethm

\brm
As before, this gives us an $H$-invariant open subvariety of $X$ with a geometric quotient by $H$ which is open in the projective variety $X \git_{\mc{L}} H$; this subvariety is the image under $\pi$ of the intersection of $\widehat{X}^{s,H}$ with the complement of the exceptional divisor in $\widehat{X}$.

In this situation the centres of the blow-ups used to obtain $\widehat{X}$ from $X$ are determined by the dimensions of the $U^{(j)}$-stabilisers  for $x \in X^{0}_{\min}$ of the limit $\lim_{t \to 0} x$. By applying Theorem \ref{Uh4} to the closures of the subvarieties where these dimensions take different values, and combining this with the partial desingularisation construction of \cite{K2} for reductive GIT quotients, $X$ can be stratified so that each stratum is a locally closed   $H$-invariant subvariety of $X$ with a categorical quotient by the action of $H$ (cf. \cite{Berczi2016a} \S 5).   This stratification can be refined further so that each stratum has a geometric quotient by the action of $H$.
\erm

\brm \label{extrm} We may use all of these theorems to perform the quotient of $X$ by $H$ using an external $\Gm$, in the following sense. We apply one of the theorems above with $H$ replaced by $\hat{H} = H \rtimes \Gm$, and $X$ replaced by $X \times \pno$, linearising the action by tensoring $\mc{L}$ with $\mc{O}_\pno(N)$ for $N \gg 0$. This yields a projective variety, which if semistability equals stability for the unipotent radical is just $(X \times \pno) \git \hat{H}$, containing an open subvariety which is a geometric quotient of a certain $H$-invariant open subset $X^{\hat{s},H} \subset X$ by $H$.   This approach is of course particularly useful if no admissible internal $\Gm$ exists, for example if $H$ is unipotent. \erm

\section{Change of distinguished $\Gm$}
 
Our investigation of the dependence of the quotient on the linearisation and one-parameter subgroup grading the unipotent radical begins with the choice of $\Gm$. For now we fix a choice of ample linearisation $\mc{L}$, and observe what happens when we vary the $\Gm$ amongst those to which $\mc{L}$ is (well) adapted. Throughout this section we assume that the condition $(*)$ holds (semistability coincides with stability for the unipotent radical) for the linearisation $\mc{L}$. For simplicity of notation we will also assume that $H$ is connected; this involves no loss of generality since GIT for finite group actions is independent of the choice of linearisation.

 The following is immediate from the Hilbert-Mumford criterion, Theorem \ref{Uh} (\ref{NRHM}).

\bpp \label{NRHM1} $X\git_\mc{L}H$ is independent of the choice of internal admissible $\Gm$ to which $\mc{L}$ is (well) adapted. \epp 

In the external case, we have no reason to suppose the quotient $(X \times \PP^1)/\!/\hat{H}$ to be independent of the choice of $\Gm$ with $\hat{H}=H\rtimes \Gm$, since we are quotienting by a group action that is dependent on that choice. However, we can use an argument  very much in the spirit of \cite{Thaddeus1996} to relate the two quotients. The plan is to reduce a change of external $\Gm$ to a change of linearisation for a certain group action on an auxiliary space.

In order to compare the quotients associated to different choices of external $\Gm$, we first wish to find an abelian group which 
 simultaneously contains every choice of admissible external $\Gm$ up to an appropriate automorphism of $X$. Then given two choices of admissible external $\Gm$ represented by one-parameter subgroups of this abelian group, we will set up an action of $(\Gm)^2$ such that the $H$-quotient with respect to each of these choices of external $\Gm$ can be viewed as a quotient of $X \times (\PP^1)^2$ by $H \rtimes (\Gm)^2$ with respect to some linearisation. We construct such an abelian group as follows. 
  
 By assumption $H$ is connected, so the action of $H$ on $X$ yields an homomorphism $\phi: H \rar \Aut^0(X)$ into the connected component $Aut^{0}(X)$ of the identity in the group $Aut(X)$ of automorphisms of $X$. Any admissible $\Gm$ must by hypothesis act on $X$ via $\psi: \Gm \to \text{Aut}^0_H(X)$, where $$\text{Aut}^0_H(X) \coloneqq \{ g \in \text{Aut}^0(X) \mid g \text{ normalises } \phi(H) \text{ and the induced action on }\phi(H)/\phi(U) \text{ is trivial} \}.$$ 
Moreover the existence of an admissible $\Gm$ implies that $\phi|_U$ is injective. 

 Now  Aut$^0(X)$ is an algebraic group (see e.g. \cite{Brion2017}), and hence so is Aut$_H^0(X)$, but Aut$_H^0(X)$ may not be linear. However, by Theorem 1 of \cite{Brion2015}, there exists a smallest normal subgroup $N_{X,H} \trianglelefteq \text{Aut}_H^0(X)$ such that the quotient $\text{Aut}_H^0(X)/N_{X,H}$ is affine, and moreoever $N_{X,H}$ is contained in the centre of $\Aut^0_H(X)$. 
  Choose a maximal torus $T_{X,H}$ of this linear algebraic group $\text{Aut}_H^0(X)/N_{X,H}$. 
  
   Given any one-parameter subgroup $\lambda : \Gm \rar \text{Aut}^0_H(X)$ of $\text{Aut}^0_H(X),$ we can conjugate by some element of $\text{Aut}_H(X)$ so that $\lambda(\Gm) \leq \pi^{-1}(T_{X,H}),$ where $\pi : \text{Aut}_H^0(X) \rar \text{Aut}_H^0(X)/N_{X,H}$ is the quotient map. Thus $\pi^{-1}(T_{X,H})$ is an abelian subgroup of $\Aut_{H}^\text{o}(X)$, such that any admissible $\Gm$ is conjugate to one in $\pi^{-1}(T_{X,H})$ via an element of $\Aut^\text{o}_H(X)$. 
   
   This means that for any admissible one-parameter subgroup $\lambda : \Gm \rar \Aut_H(X)$, there is an automorphism $\a :X \rar X$ conjugating $\lambda$ to  $\m : \Gm \rar \pi^{-1}(T_{X,H})$, along which any linearisation $\mc{L}$ of the $\Uh_\lambda$ action may be pulled back, so that $\a$ induces an isomorphism $$X \git_\mc{L} \Uh_{\lambda} \cong X \git_{\a^*\mc{L}} \Uh_{\m}.$$ Finally, we observe that no admissible $\Gm$ can lie in $N_{X,H}$: since the latter is contained in the centre of $\Aut^0_H(X)$, any $\Gm$ it contains would commute with the $H$-action, and hence have trivial conjugation action on $U$. This means that any admissible $\Gm$ may, up to automorphism of the whole picture, be considered to be a subgroup of $T_{X,H}$.

We can now investigate the dependence on the choice of external $\Gm$. Consider the action of $\Hh$ on $X \times \pno$, linearised so as to weight the $\pno$ factor heavily, as  in  Remark \ref{extrm}. We will show that, under certain circumstances, changing the one-parameter group $\Gm$ for the action of $H\rtimes \Gm$ on $X$ may equivalently be viewed as changing the linearisation for an action $$(H \rtimes (\Gm \times \Gm)) \times (X \times \pno \times \pno) \rar X \times \pno \times \pno.$$ \\ We begin with a linearisation $\mc{L}^{\prime}$ of the $H$ action on $X$, with respect to an underlying line bundle $L^{\prime}$. Let $T_{\lambda}$ and $T_{\mu}$ be two copies of $\Gm$ equipped with actions on $X$ and $H$, admissible respectively for linearisations $\mc{L}_\lambda $ and $\mc{L}_{\mu}$ of the actions of $\hat{H}_\lambda \coloneqq H \rtimes T_{\lambda}$ and $\hat{H}_\mu \coloneqq H \rtimes T_{\mu}$ on $X \times \pno$, with respect to the same underlying bundle $L = L^{\prime}\boxtimes \mc{O}_{\pno}(N)$, so that the induced linearisations for the $H$-action are the same.  By the above we may assume that both $\Gm$'s lie in $T_{X,H}$ and hence commute, so that we get an induced action of their product on $X$.  The goal is to produce an action of $T_\lambda \times T_\mu$ on $\pno \times \pno$ and two linearisations $\mc{\hat{L}}_\lambda$ and $\mc{\hat{L}}_\mu$ for the resultant action of $H \rtimes (T_\lambda \times T_\mu)$ on $X \times \pno \times \pno$, such that we have 
\begin{equation} (X \times \pno \times \pno) \git_{\mc{\hat{L}}_{\lambda}} (H \rtimes (T_\lambda \times T_\mu)) \cong (X \times \pno) \git_{\mc{L}_{\lambda}} \hat{H}_{\lambda},  \label{1}
  \end{equation}

\begin{equation} (X \times \pno \times \pno) \git_{\mc{\hat{L}}_{\mu}} (H \rtimes (T_\lambda \times T_\mu)) \cong (X \times \pno) \git_{\mc{L}_{\mu}} \hat{H}_{\mu}. \label{2}
 \end{equation}
Having done this, it will only remain to explore the dependence of the quotient on the choice of linearisation.

\bpp \label{pp1} Adopt the notation above. Then there are linearisations $\mc{\hat{L}}_{\lambda}$, $\mc{\hat{L}}_{\mu}$ of the action $$(H \times T_{\lambda} \times T_{\mu}) \times (X \times \pno \times \pno) \rar (X \times \pno \times \pno)$$ which extend the respective linearisations on $X \times \pno$ and give isomorphisms  (\ref{1}) and (\ref{2}). \epp

\bpf Extend the actions of $H\rtimes T_{\lambda}$ and $H\rtimes T_{\m}$ to actions on $X\rtimes \pno_{\lambda} \rtimes \pno_{\mu}$ by declaring that $T_\lambda$ acts trivially on $\pno_{\mu}$ and vice-versa. We thus obtain two linearisations of the action of $H\rtimes(T_{\lambda} \times T_{\mu})$ with repect to the line bundle $$\hat{L} \coloneqq L \boxtimes \mc{O}_{\pno_{\lambda}}(N)\boxtimes \mc{O}_{\pno_{\mu}}(N)$$ on $X \times \pno_{\lambda} \times \pno_{\mu}$: one by tensoring $\mc{L}_\lambda$ with the pullback of the trivial linearisation on $\mc{O}_{\pno_{\mu}}$, and the other with the roles of $\lambda$ and $\mu$ reversed. We tensor these together, and call the result $\widehat{\mc{L}}$.

Now consider the weight lattice of $T_\lambda \times T_\mu$ with respect to this linearisation. Using a Segre embedding  into a  projective space $X \times \pno_{\lambda} \times \pno_{\m} \hrar \mathbb{P}^s$, we find that for the linearisation decribed above the weights are arranged in four clusters, each placed at one of the four corners of a square. Each cluster is a small copy of the weights of $X$, and the weights in it correspond to taking the coordinates of $X$ in conjunction with one of the four pairs $y_iz_j$ for $i,j=0,1$, where $(y,z)\in \pno_{\lambda} \times \pno_{\m}$. Let $\widehat{\mc{L}}_{\lambda}$ be  the twist of $\widehat{\mc{L}}$ by the $T_\lambda \times T_{\m}$-character $(0,N+r_{\lambda}-\eps)$, and let $\widehat{\mc{L}}_{\lambda}$ be  the twist of $\widehat{\mc{L}}$ by the character $(N+r_{\m}-\eps,0)$, where $r_{\lambda}$ and $r_{\m}$ are respectively the minimal weights on $\mc{L}_{\lambda}$  and $\mc{L}_{\m}$ for their distinguished choices of $\Gm$.

Consider the quotient with respect to $\hat{\mc{L}}_{\lambda}$; the quotient with respect to $\hat{\mc{L}}_{\mu}$ works in the same way. By Proposition \ref{NRHM1}, the internal quotient is independent of the choice of 1PS used to form $\Uh$, so we may choose  $\Uh$ to be $\Uh_{\lambda} = U \rtimes T_{\lambda}$. Then $\Uh_{\lambda}$ acts trivially on the third factor, and by our choice of character, the quotient is $$(X \times \pno_{\lambda} \times \pno_{\m}) \git_{\hat{\mc{L}}_{\lambda}} \Uh_{\lambda} \cong ((X\times \pno_{\lambda})\git_{\mc{L}_{\lambda}} \Uh_{\lambda}) \times \pno_{\m}.$$ We follow this by the residual reductive quotient by $R \times T_{\mu}$, which we may factor as first a quotient by $T_{\mu}$, then a quotient by $R$. We observe that the $T_{\mu}$ quotient is just the first factor of the product, since the stable (and semistable) locus is $$(((X\times \pno_{\lambda})\git_{\mc{L}_{\lambda}} \Uh_{\lambda}) \times \pno_{\m})^{s,T_{\m}} = ((X\times \pno_{\lambda})\git_{\mc{L}_{\lambda}} \Uh_{\lambda}) \times (\mathbb{A}_{\m}^1 \setminus \{0\}).$$ In addition, the induced linearisation of the $R$-action on this first factor is equal to (a positive scalar multiplie of) the $R$-linearisation induced by $\mc{L}_{\lambda}.$ So we have $$((X\git_{\mc{L}_{\lambda}} \Uh_{\lambda}) \times (\pno_{\mu}))\git (R \times T_{\m}) \cong (X\git_{\mc{L}_{\lambda}} \Uh_{\lambda}) \git R \cong (X \times \pno) \git_{\mc{L}_{\lambda}} \hat{H}_{\lambda}.$$ 
This gives us the isomorphism  (\ref{1}); the isomorphism (\ref{2}) follows in the same way.

\epf

\bcr  Changing between two admissible external $\Gm$'s, for a fixed linearisation well-adapted to both, is equivalent to a change of linearisation. \ecr

\section{Change of linearisation}

Continuing with our analysis, we now tackle the other half of the question. We continue to assume that we have semistability coincides with stability for the unipotent radical. Fix an admissible one-parameter group $\Gm$, internal or external, and consider the space of linearisations well-adapted to it. The following observation is immediate from Theorem  \ref{Uh}(3).

\bpp \label{prop41}
Let $\mc{L}_1$ and $\mc{L}_2$ be two well-adapted linearisations of the $H$-action on $X$, and suppose that their restrictions to a maximal torus $T$ of $H$  give the same $T$-semistable locus. Then $X^{ss,H}(\mc{L}_1) = X^{ss,H}(\mc{L}_2)$, and the quotients $X\git_{\mc{L}_1} H$ and $X\git_{\mc{L}_2}H$ can be canonically identified.
\epp

This gives us part of the story: if we consider the structure of the space of linearisations for the action of the maximal torus $T$ of $H$, and identify a linearisation of the $H$-action with the linearisation of the $T$-action that it induces, then by Proposition \ref{prop41}, if two linearisations lie in the same chamber (or more generally are GIT-equivalent for the action of $T$ in the sense of Dolgachev \& Hu), then they must produce the same $H$-quotient. What we have not yet determined is how the quotient may change when we cross a wall (or more generally move from one GIT-equivalence class to another). We consider first the situation when $H=\Uh$.

\bpp The quotient $X\git_{\mc{L_i}}\Uh$ is independent of the choice of ample well-adapted linearisation $\mc{L_i}$.
\label{indepUh}

\epp

\bpf Recall that the $\Uh$ quotient admits a description as the geometric quotient of $$X^{s,\Uh} = X \setminus U(X \setminus X^{s,\Gm})$$ by $\Uh$, where $X^{s,\Uh}$ is an open subvariety of $X$ dependent only on the $\Uh$-action and the stable locus for our chosen one parameter group $\lambda:\Gm \to \Uh$ which grades $U$, and any two choices for $\lambda$ are conjugate in $\Uh$. It suffices to show that $X^{s,\Gm}$ is independent of the choice of linearisation $\mc{L}$.

 One way to argue this is the following (compare Theorem 4.3 of \cite{Bialynicki-Birul1973} on the Bia\l ynicki-Birula stratification). The one-parameter subgroup $\lambda :\Gm \rar \Uh$ determines a map of sets $p : X \rar X$ defined by $$p(x) \coloneqq \lim_{t \rar 0}\lambda(t)x \in X$$ taking each point to its limit under the flow given by the action of $\lambda:\Gm \to \Uh$, and $p(X) = X^{\Gm}$ is the fixed point set for the action of $\lambda$. 
 
Let $Z_1, \ldots,Z_n$ be the connected components of $X^{\Gm}$. By assumption $X$ is irreducible, so $X_{\min}$ is the union of those $Z_j$ for which $p^{-1}(Z_j)$ is open in $X$, and $X^{s,\Gm} = p^{-1}(X_{\min}) \setminus X_{\min}$ is independent of $\mc{L}$. 
\epf

In general,  two different $H$-linearisations will induce different linearisations of the residual $\bar{R} \coloneqq R/\Gm$ action on $X \git \Uh$, and hence the two $H$ quotients for change of linearisation will be related by a reductive VGIT picture applied to $X \git\Uh$. This allows us to recover, for fixed $\Gm$, a wall-and-chamber structure to the space of well-adapted linearisations, namely that of the reductive VGIT applied to the $T$-action, so that two quotients are related by Thaddeus flips of $X\git \Uh$.

\section{When semistability does not coincide with stability}
\label{ssns}
So far we have been considering the situation covered by Theorems \ref{Uh} and \ref{Uh3} when we can construct quotients of the variety $X$ (or $X\times \pno$ in the external case) without needing to do the blow-up process of Theorems \ref{Uh2} and \ref{Uh4}. As  described in \cite{Berczi2016a}, this process begins by defining $$d_0 = \max \{\dim \stab_U(x) \mid x \in X^{0,\Gm}_{\min}\}, $$ and then blowing up along the closure in $X$ of $$\Delta_0 = \{ x  \in X^{0,\Gm}_{\min} \mid \dim \stab_U(x) = d_{0} \}. $$ Iteratively, one then sets $X_{(0)} = X$ and defines $X_{(j)}$ to be the blow-up along the closure in $X_{(j-1)}$ of $$\Delta_j = \{ x  \in X^{0,\Gm}_{\min} \mid \dim \stab_U(x) = d_j \},$$ where $$d_j = \max \{\dim \stab_U(x) \mid x \in (X_j)^{0,\Gm}_{\min}\}. $$ This procedure is clearly dependent on the choice of $\Gm$, as written. However, one may describe it in a different way so that the results above will apply to these sorts of quotients, too. Instead of blowing up along the loci described, we let $$X^0 = \bigcup_{\Gm \in \mathfrak{C}} X^{0,\Gm}_{\min}. $$ This is the set of points in $X$ that are in the basin of attraction to the minimal weight space for some admissible $\Gm$. Then we can apply the procedure described in \cite{Berczi2016a} with $X^{0,\Gm}_{\min}$ replaced by $X^0$ to achieve a \lq universal\rq\text{ }blow-up, $\check{X}$, that is not dependent on any choice. For each admissible one-parameter group $\Gm$ in $\mathfrak{C}$, there is an open subvariety of $\check{X}$ that is isomorphic to the open subset $(\hat{X}_{\Gm})^{0,\Gm}_{\min} \subset \hat{X}_{\Gm}$, where the latter is the space obtained from $X$ by applying the procedure above using $X^{0,\Gm}_{\min}$. All quotients obtained from the blow-up process using different choices of $\Gm$ can be regarded as quotients of the \lq master space' $\check{X}$ using different choices of $\Gm$, and hence the quotients are related as already described in the cases covered by Theorems \ref{Uh} and \ref{Uh3}.

\section{Summary of results}

As we have seen, the main difference between the parameter spaces for reductive and non-reductive GIT is the additional choice of a one-parameter group $\Gm$. This can be chosen from the admissible cone $\mathfrak{C} \subset \mathfrak{z}$ in the internal case, or from a corresponding cone $\mathfrak{C} = \mathfrak{C}_{X,H} $ in the external case. Without loss of generality we may assume that semistability coincides with stability for the unipotent radical; if not we can work with the universal blow-up described in $\S$ \ref{ssns}. Having made a choice of one-parameter group $\Gm$, we must choose a linearisation well-adapted to this choice.

\bdf Given a choice of admissible internal $\lambda :\Gm \rar H$, let $\partial_{\lambda} \subset \mathfrak{C}$ be the set of linearisations that are well-adapted to $\lambda$. Let $$ \partial = \bigcup_{\lambda \mbox{ admissible}} \partial_\lambda
$$ be the set of linearisations that are well-adapted for some admissible internal $\Gm$. \edf

The set $\partial_{\lambda}$ is a suitably small open region along the boundary of the intersection of $\mathfrak{C}$ with a region bounded by several hyperplanes: one giving the negativity condition for the minimal weight space for $\lambda: \Gm \to H$ (which by Proposition \ref{indepUh} is the same for all ample linearisations), and the others giving the positivity conditions for the other weight spaces. In particular, $\partial_{\lambda}$ is the intersection with a convex set of an open neighbourhood of its boundary in $\mathfrak{C}$. While this region evidently depends on the choice of $\Gm$, we can recover $\partial$ as the union of $\partial_{\lambda}$ for only finitely many $\lambda$. 

We may summarise our results as follows.

\bthm \textbf{(Internal Case)}
Let $X$ be a projective variety, acted on by a linear algebraic group $H$. 
\bnu
\item Suppose that condition $(*)$ or $(***)$ is satisfied. Then the quotient $X\git_{\mc{L}}H$ is independent of the choice of admissible internal $\Gm$, and may be formed using any linearisation $\mc{L} \in \partial$ which is well-adapted to some admissible $\Gm$. Moreover if two linearisations are GIT-equivalent for the action of a maximal torus $T \leq H$ then the associated quotients are isomorphic. Further, any two quotients with respect to different  linearisations, which are both well-adapted for the same one-parameter subgroup, are related by the usual reductive variation of GIT picture induced via the $\Uh$ quotient map for the relevant actions of $R/\Gm$ on the quotient $X \git \Uh$ (where $X \git \Uh$ is independent of the choice of linearisation).
\item If we do not have $(*)$ or $(***)$, then there is a sequence of $H$-equivariant blow ups of $X$ along closed $H$-invariant subvarieties resulting in a variety $\widetilde{X}$ such that either $(*)$ or $(***)$ holds, as desired, and the results of the first two parts of this Theorem apply.
\enu
\ethm

Note in particular that there are only finitely many possible GIT quotients yielded by the $\Uh$ theorems with internal $\Gm$, and we may move between any two of these by a sequence of changes of linearisation, if necessary changing the $\Gm$ to maintain well-adaptedness (though without changing the quotient). We therefore interpolate between any two $H$-quotients by a sequence of \lq vertical\rq\text{ }and \lq horizontal\rq\text{ }moves in $\mathfrak{C}\times \partial$, giving a sequence of points such that any two consecutive points have either the linearisation or the one-parameter subgroup $\Gm$ in common.

\bthm \textbf{(External case)} 
Let $X$ be a projective variety, acted on by a linear algebraic group $H$. Form the non-reductive GIT quotient 
$$ X \hat{\git} H = (X \times \PP^1) \git \hat{H}$$
with respect to a choice of admissible external $\Gm \leq \Aut H$ and a well-adapted ample linearisation $\mc{L}$ by letting $\hat{H} = H\rtimes \Gm$ act on $X\times \pno$, linearised with respect to $\mc{L}\boxtimes \mc{O}_{\pno}(N)$ for $N \gg 0$. 
\bnu
\item There is an abelian subgroup  $\pi^{-1}(T_{X,H})$ of $\text{Aut}(X)$ such that any admissible external $\Gm$ is GIT-equivalent to a one-parameter subgroup $\Gm$ of $\pi^{-1}(T_{X,H})$. Furthermore external one-parameter groups $\Gm$ yield the same non-reductive quotient $ X \hat{\git} H$ for the same linearisation if they have the same minimal weight space in $X$. 
\item The dependence of the quotient $ X \hat{\git} H$ on the linearisation can be described as in the previous theorem.
\enu

\ethm

\section{Examples and Applications}

We conclude with a brief description of some applications, beginning with a concrete example to illustrate these ideas in practice.


\subsection{Points and lines in $\mathbb{P}^2$}

We consider configurations of points and lines in the projective plane up to an action of a certain non-reductive subgroup of SL$_3(\CC)$. We take $V$ to be a 2-dimensional complex vector space, and let $$X= \mathbb{P}(V)^p \times \mathbb{P}(V^*)^q,$$ so that a point of $X$ corresponds to a choice of $p$ points and $q$ lines in the projective plane. Our first choice for the group in question is $$H = \Bigg\{ \begin{pmatrix}
a & b & c \\ 0 & d & 0 \\ 0 & 0 & f
\end{pmatrix} \in \text{SL}_3(\CC) \Bigg\},$$ and the unipotent radical and our chosen maximal torus are, respectively,
$$U = \Bigg\{ \begin{pmatrix}
1 & b & c \\ 0 & 1 & 0 \\ 0 & 0 & 1
\end{pmatrix} \in \text{SL}_3(\CC) \Bigg\}$$

$$T = \Bigg\{ \begin{pmatrix}
t & 0 & 0 \\ 0 & s & 0 \\ 0 & 0 & (st)^{-1}
\end{pmatrix} \in \text{SL}_3(\CC) \Bigg\}.$$ Note that in this example $H$ is just the semidirect product of a torus and the unipotent part, so the torus is a quotient group and hence all its characters extend to the whole group $H$. The action of $H$ on $X$ is that induced by the natural action on $V$, and the dual action on $V^*$. Accordingly, the weight polytope is scaled version of the hexagon familiar from the study of the representation theory of SL$_3$. Take $p=2$, $q=1$ for concreteness.

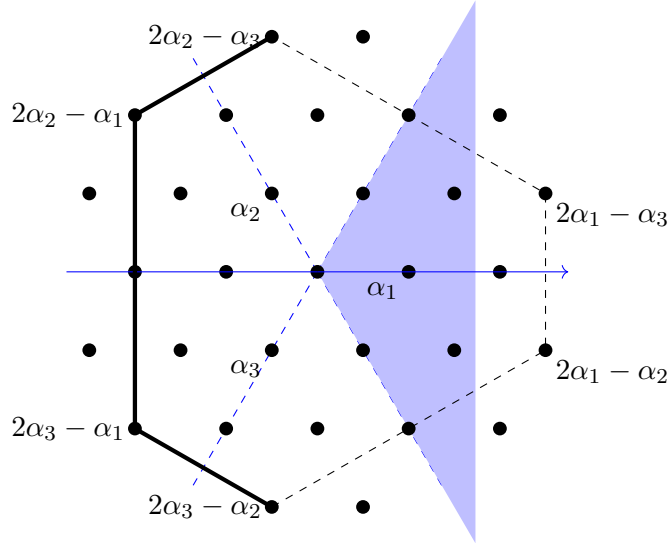
\begin{figure}[H] \label{weightsone}
\centering
\caption{Weight diagram and admissible cone for the action of $H$ on $(\PP^2)^2 \times (\PP^2)^*$}
\begin{tikzpicture}[scale=0.6]

\pgfmathsetmacro\ax{2}
\pgfmathsetmacro\ay{0}
\pgfmathsetmacro\bx{2 * cos(120)}
\pgfmathsetmacro\by{2 * sin(120)}
\pgfmathsetmacro\lax{2*\ax/3 + \bx/3}
\pgfmathsetmacro\lay{2*\ay/3 + \by/3}
\pgfmathsetmacro\lbx{\ax/3 + 2*\bx/3}
\pgfmathsetmacro\lby{\ay/3 + 2*\by/3}

\foreach \k in {1,2,4,5} {
  \draw[blue,dashed] (0,0) -- +(\k * 60:5.5);
}
\fill[blue!25] (0,0) -- (-60:6 * \lby) -- (60:6 * \lby) -- cycle;
\begin{scope}
\clip (90:6) \foreach \k in {1,...,6} { -- ++([rotate=\k * 60 + 60]90:6) };
\foreach \na in {-3,...,3} {
  \foreach \nb in {-3,...,3} {\
    \node[circle,fill=black,scale=0.5] at (\na * \ax + \nb * \bx, \na * \ay + \nb * \by) {};
  }
}
\end{scope}
\draw[blue] (0,0)--(3*60:5.5);
\draw[blue, ->] (0,0)--(6*60:5.5);
\node at (\ax,\ay) [below left] {\(\alpha_1\)};
\node at (\bx,\by) [below left] {\(\alpha_2\)};
\node at (\bx,-\by) [below left] {\(\alpha_3\)};
\node at (2.5*\ax,\by) [below right] {\(2\alpha_1-\a_3\)};
\node at (2.5*\ax,-\by) [below right] {\(2\alpha_1-\a_2\)};
\node at (\bx,3*\by) [left] {\(2\alpha_2-\a_3\)};
\node at (-2*\ax,2*\by) [left] {\(2\alpha_2-\a_1\)};
\node at (-2*\ax,-2*\by) [left] {\(2\alpha_3-\a_1\)};
\node at (\bx,-3*\by) [left] {\(2\alpha_3-\a_2\)};
\draw[black, dashed](\bx,-3*\by)--(2.5*\ax,-\by)-- (2.5*\ax,-\by) -- (2.5*\ax,\by) -- (\bx,3*\by);
\draw[black, ultra thick] (\bx,3*\by)--(-2*\ax,2*\by)--(-2*\ax,-2*\by)--(\bx,-3*\by);

\end{tikzpicture}
\end{figure}

The admissible cone, shaded in blue, is larger than the standard positive Weyl chamber in $\SL (3)$, because we do not have the full Borel subgroup, and so we only need positivity with two roots (cf. \cite{Berczi2016a} 7.1). The well-adapted linearisations correspond to the section of the boundary indicated in black.  We observe also that for this $H$ there is a natural choice for the 1PS, given by the solid blue arrow: a universal choice in the sense that any linearisation which is allowable for some admissible $\Gm$ will be allowable for this one. It should be noted that this does not hold in general. For example, if instead we impose also the additional condition that $b=0$, so that $H$ is the semidirect product of $\Ga$ and the maximal torus, the admissible cone is enlarged and there is no longer a universal choice of 1PS.
%
%
%
%
%

\subsection{Unstable strata}
A final source of applications is the \lq unstable strata\rq\text{ }problem. For this recall the definition of the Hilbert--Mumford function.

\bdf Let $X \subset \PP^r$ be a projective scheme, $G$ a reductive linear algebraic group acting on $X$. Choosing a linearisation of this action with respect to a very ample line bundle, we may assume the action is via a representation $\phi : G \rar \GL_{r+1}$. Take a 1PS $\rho : \Gm \rar G$, choose coordinates to diagonalise the action of $\rho$, and let $w=(w_0,..,w_r)$ be its weight vector. For a point $x = [x_0:..:x_r] \in X$  we define the Hilbert-Mumford function to be $$\mu (x,\rho) = \max \{ -w_i \mid x_i \neq 0 \}.$$ \edf

For our purposes it will be convenient to work with a normalised version. 

\bdf Fix an invariant inner product on the set of 1PS's of $G$: that is, take a positive definite Weyl-invariant integer-valued bilinear form on the set of 1PS's of any maximal torus $T$. Denote by $\|-\|$ the associated norm. Then define the normalised Hilbert-Mumford function to be
 $$M(x,\rho) = \frac{\mu(x,\rho)}{\| \rho \|}.$$ 
 \edf

\brm If the group $G$ is semisimple, for example if $G = \SL_{r+1}$, we have a natural choice of invariant inner product, i.e. the Killing form. \erm

\bdf We say that a 1PS $\rho : \Gm \rar G$ is a \emph{maximally destabilising one-parameter subgroup} (md1PS) of a point $x \in X$, or that $\rho$ is \emph{adapted} to $x$, if $M(x,\rho)$ is minimal amongst all 1PS's of $G$. \edf 

\bex 
If $G=T$ is a torus diagonalised by coordinates $x_0,..,x_r$ with weights $\a_0,..,\a_r$, then the minimal value of the normalised Hilbert-Mumford function for $x = [x_0:..:x_r]$ is the norm of the closest point to the origin, $B(x)$, of the convex hull of the $T$-weights in the state polytope of $x$, i.e. those $w_i$ such that $x_i \neq 0$. There is, up to scaling, a unique maximally destabilising 1PS for $x$, namely that obtained by exponentiating $B(x)$.
\eex

The following result tells us what to expect more generally.

\bpp \label{adapt} \cite{Kempf1978} Let $G$ be a reductive group scheme acting on a projective scheme $X$ linearly with ample linearisation $\mc{L}$. Let $x \in X$ be an unstable point. Then
\bnu 
\item There exists $\rho : \Gm \rar G$ which minimises the function $M(x,-)$ amongst all $1$PS's of $G$.
\item There exists a parabolic subgroup $P_x$ of $G$ such that for any $\lambda$ satisfying $(1)$ we have $$P_x = P(\lambda) \coloneqq \{ g \in G  \mid \lim_{t \rar 0} \lambda (t) g \lambda (t^{-1}) \text{ exists in $G$} \}. $$
\item Any $1$PS satisfying $(1)$ is unique up to replacing $\lambda$ by $\lambda^n$ and conjugation by an element of the subgroup $P_x $.
\enu
\epp
 
This enables us to associate to an unstable point a 1PS that is \lq most responsible\rq\text{ }for that point's instability, unique in the sense described above. In applications to moduli problems, these are the subgroups that we think of as somehow best highlighting whatever intrinsic property of the object in question makes it moduli-unstable. \\~\\ We now define the $\beta$-stratification associated to a linear action of a reductive linear algebraic group on a projective variety.

\bdf Fix a maximal torus $T$ of $G$ acting as above, and a positive Weyl chamber $\mathfrak{t}^+$ in its Lie algebra. Define the index set for the stratification, $$\mc{B} \coloneqq \{\beta \in \mathfrak{t}^+ \mid \exists x \in \mathbb{P}^r \text{ s.t. }B(x) = \beta \}.$$ \edf To put it another way, $\mc{B}$ is the set of points in the positive Weyl chamber that are the closest point to the origin of the convex hull of some non-empty subset of the $T$-weights. To each element of this index set we associate a 1PS.

\bdf For an element $\beta \in \mc{B}$, denote by $\lambda_{\beta}$ the 1PS that is given by the weight vector $q\beta$, where $q$ is the smallest positive rational number giving integral weights. \edf

Now that we have all the combinatorial information we will need, we can begin to show how this data captures the instability of all points of $X$. Given $\beta \in \mc{B}\setminus \{0\}$, let $$H_{\beta} = \{v \in \mathfrak{t} \mid v\cdot\beta = \|\beta\|^2 \},$$ $$H_\beta^{+} = \{ v \in \mathfrak{t} \mid v \cdot \beta \geq \|\beta\|^2 \}$$ denote respectively the hyperplanes through $\beta$ meeting it perpendicularly, and the half-space on the opposite side of this hyperplane from the origin.

\bdf Following \cite{Kirwan1984}, we define $$\Zb \coloneqq \{ x \in X \mid \text{ all weights of $x$ lie on } H_{\beta} \} $$ $$\Yb \coloneqq \{ x \in X \mid \text{ all weights of $x$ lie in $H^+_{\beta}$ and at least one weight lies on }H_{\beta} \}$$ \edf

The first of these is a closed subvariety of $X$, and the second is a locally closed subvariety. Clearly both contain only unstable points for all $\beta \in \mc{B} \setminus \{0\}$. There is a natural surjective retraction $p_\beta : \Yb \rar \Zb$, given by taking the limit as $t \rar 0$ under $\lambda_{\beta}$. If $X$ is nonsingular the fibres are just affine spaces, by \cite{Bialynicki-Birul1973}.
\bdf
We further define $$\Yb^{ss} = \{ x \in \Yb \mid \lambda_{\beta} \text{ is adapted to } x \}, $$  $$\Zb^{ss}  = \{x \in \Zb \mid \lambda_{\beta} \text{ is adapted to } x \},$$ $$S_{\beta} = G\cdot \Yb^{ss}.$$
\edf  Denote by $\stab\beta$ the stabiliser of $\beta$ under the adjoint action of $G$, and note that $\Zb$ is invariant under the action of this group. We restrict the linearisation to the action of $\stab \beta$ on $\Zb$, and twist it by the character $-\beta$; this has the effect of shifting the weights by this vector. Then $\Zb^{ss}$ is exactly the semistable locus for the action of $\stab \beta$ on $\Zb$ with respect to the twisted linearisation, and $\Yb^{ss} = p_\beta^{-1}(\Zb^{ss})$.

\bdf
Given $\beta \in \mc{B}$, denote by $P_{\beta}$ the parabolic subgroup of $G$ associated to $\lambda_{\beta}$. Denote by $S_{\beta}\times^{P_{\beta}}G$ the quotient of $S_{\beta}\times G$ by the free action of $\Pb$ given by $h\cdot(g,x) \coloneqq (gh^{-1},hx).$
\edf

\brm Though a priori an algebraic space, by \cite{Popov1994} we know that $S_{\beta}\times^{P_{\beta}}G$ is in fact an algebraic variety. \erm

We can now summarise the relevant results of \cite{Kirwan1984} that give the stratification.

\bthm \cite{Kirwan1984}
Let $X$ be a projective variety, acted on by a reductive linear algebraic group $G$, linearly with respect to an ample line bundle $\mc{L}$. Then, with notation as above, we have 
\bnu
\item There is a stratification of $X$ into disjoint locally closed $G$-invariant subvarieties, $$X = \bigsqcup_{\beta \in \mc{B}} S_{\beta},$$  where the open stratum $S_0 = X^{ss}$, and the closure of a stratum $S_{\beta}$ satisfies $$\overbar{S}_{\beta} \subset \bigcup_{\substack{\gamma \text{ s.t.} \\ \|\gamma\| \geq \|\beta\|}} S_{\gamma}.$$
\item There is, for every $\beta \in \mc{B}$, a $G$-equivariant isomorphism $S_{\beta}\times^{P_{\beta}}G \cong \Yb^{ss}.$
\enu
\ethm 

We can try to find, for every $\beta \in \mc{B}$, a quotient of (an open subvariety of) the stratum $S_{\beta}$ by the $G$-action. By ($2$) in the theorem above, this is equivalent to the problem of finding a quotient for the action of $P_{\beta}$ on (an open subvariety of) $\Yb^{ss}$. 

There is in fact a categorical quotient for the  action of $P_{\beta}$ on $\Yb^{ss}$ given by the reductive GIT quotient of $\Zb$ by $\stab\beta$, where the linearisation is twisted by $\beta$. However this categorical quotient is very far from being a geometric quotient. 
 If we are to obtain quotients of the correct dimension, we must prevent points of $\Yb^{ss}$ from being identified with their images under $p_{\beta}$ by making $\Zb^{ss}$ unstable. We therefore split each stratum into two pieces: roughly speaking, one piece corresponds to $\Zb^{ss}$, and one to the open subset of $\Yb^{ss}$ not in the $P_{\beta}$ sweep of $Z^{ss}_{\beta}$. The GIT quotient of $\Zb$ by $\stab\beta$ is a categorical  quotient of the former, and we can using Theorem \ref{Uh} and its variants to obtain quotients of open subvarieties of the latter.
 \\~\\ 
In particular, given a GIT-formulation of a moduli problem, we can try to interpret the $\beta$-stratification in terms of the intrinsic properties of the objects themselves, and form quotients of each $\beta$-stratum to get moduli spaces of objects of fixed $\beta$-type.
This applies to some classical moduli problems with well-known GIT formulations: moduli of semistable sheaves (or complexes of sheaves) over a projective base, and moduli of stable curves. In each case, one must first understand the $\beta$-stratification, associating to each sheaf or curve its \emph{$\beta$-type}, which says in what level of the stratification it may be found\footnote{Of course, these objects come in unbounded families, and in any given GIT-setup we will only see finitely many strata, so strictly speaking we must work with the asymptotic $\beta$-stratification. See \cite{Berczi2017} for more details.};  for sheaves, this involves the Harder-Narasimhan type, which records the Hilbert polynomials of the subquotients in the Harder-Narasimhan filtration (cf. \cite{Hoskins2012}). This gives a discrete classification, after which the second task is to apply $\S$2 to construct moduli spaces of unstable objects with fixed discrete invariants.

\bibliography{con2}
\bibliographystyle{alpha}

\end{document}